\documentclass[12pt]{article}
\usepackage{amssymb}
\usepackage{amsfonts}

\input{tcilatex}
\begin{document}

\begin{center}
\bigskip

\textbf{Note on antichain cutsets in discrete semimodular lattices}
\end{center}

\bigskip

\begin{center}
Stephan Foldes

Tampere University of Technology

PL 553, 33101 Tampere, Finland

sf@tut.fi

February 2011

\bigskip

\textbf{Abstract}
\end{center}

\textit{The characterization of level sets of finite Boolean lattices as
antichain cutsets, due to Rival and Zaguia, is seen to hold in all discrete
semimodular lattices.}

\textit{\bigskip }


\bigskip

\bigskip

\textbf{1 Background}

\bigskip

An \textit{antichain cutset} in a partially ordered set is a set of elements
intersecting every maximal chain in a singleton. Finite non-empty posets
always have antichain cutsets. In \textit{ranked posets with a least element}
$0$ (i.e. where in every interval $[0,x]$ all maximal chains have the same
finite number of elements $\neq 0$, called the \textit{height} $h(x)$ of $x$%
), each set

\[
L_{n}=\left\{ x:h(x)=n\right\} \text{ \ \ \ \ }n=0,1,...
\]%
is an antichain cutset if $L_{n}\neq \emptyset .$ An early study involving
antichain cutsets is by Grillet [G]. Rival and Zaguia have shown in [RZ],
Theorem 4, that in finite Boolean lattices height classes $L_{n}$ are the
only antichain cutsets. In the next section we shall see that this result
extends to all discrete semimodular lattices. Such lattices may not have a
least element, we therefore define, in any poset, a \textit{level} as an
equivalence class of the equivalence relation $\equiv $ which is obtained as
the reflexive-transitive closure of the following symmetric relation $\sim $
:%
\[
x\sim y\text{ \ \ }\Leftrightarrow \text{ \ \ }\exists \text{ }z\text{
covered by both }x\text{ and }y
\]%
For finite Boolean lattices this gives an alternative description of the
sets $L_{n}$\ which remains meaningful in the larger context of discrete,
possibly infinite and unbounded semimodular lattices. By \textit{%
semimodularity} we understand the lower covering condition 
\[
x\text{ covers }x\wedge y\text{ \ }\Rightarrow \text{ \ }x\vee y\text{
covers }y
\]%
By a \textit{discrete} order we mean a poset in which every interval $[x,y]$
has a finite maximal chain. In a discrete semimodular lattice, for $x\leq y$
all maximal chains of $[x,y]$ have the same finite number of elements $\neq x
$, called the \textit{height of} $y$ \textit{above} $x$, denoted $h(x,y)$. 

It is easy to see that in finite semimodular lattices, which are ranked
posets, levels defined as equivalence classes of the reflexive-transitive
closure $\equiv $ of the relation $\sim $ coincide with the non-empty sets $%
L_{n}=\left\{ x:h(x)=n\right\} $, \ $n=0,1,...$ In all discrete semimodular
lattices, two elements $x,y$ are in the same level class if and only if they
have the same height above their meet. In fact for any common lower bound $z$
of $x$ an $y$, $h(z,x)=h(z,y)$ if and only if $x$ and $y$ are in the same
level class.

\bigskip

\bigskip

2 \textbf{Generalized statement and proof}

\bigskip

The following generalizes Theorem 4 of [RZ].

\bigskip

\textbf{Theorem}\textit{\ \ Let }$L$ \textit{be any discrete semimodular
lattice, and let }$A\subseteq L$. \textit{Then }$A$\textit{\ is an antichain
cutset if and only if }$A$\textit{\ is a level class. }

\bigskip

\textbf{Proof \ }For $x<y$ we write $h(y,x)$ for the negative of the height
of $y$ above $x$. Thus $h(y,x)=-h(x,y)$ for all comparable $x,y.$

\bigskip

Suppose $A$ is a level class. First, if $x,y\in A$, then $h(x\wedge
,x)=h(x\wedge y,y),$ which rules out $x<y$. Thus $A$ is an antichain.
Second, let $C$ be a maximal chain. Choose any $a\in A$ and $y\in C$, and
let $z=a\wedge y$. The chain $C$ must contain an element $x$ such that $%
h(x,y)=h(z,y)-h(z,a).$ Then $a$ and $x$ have the same height above $z$ and
therefore $x$ is also in the level class $A$. This concludes the proof that $%
A$ is an antichain cutset.

\bigskip

Suppose that $A$ is an antichain cutset but not a level class. Choose any $%
a\in A$ and let $N$ be the level class of $a$. As $N$ is an antichain
cutset, $N\nsubseteq A$. Choose any $b\in A\setminus N.$ There is a sequence
of elements of $N,$%
\[
x_{0}=a,\text{ }x_{1},...,x_{n}=b 
\]%
such that $x_{i}\wedge x_{i+1}$ is covered by $x_{i}$ and $x_{i+1}$ for each 
$i=0,...,n-1.$ For the first index $i$ with $x_{i+1}\notin A$, write 
\[
x=x_{i},\text{ }y=x_{i+1},\text{ \ }z=x\wedge y\text{, \ }w=x\vee y 
\]%
As $z$ is covered by $x$ and $y$, by semimodularity $w$ covers both $x$ and $%
y$. Clearly the chain $\left\{ z,y,w\right\} $ avoids $A.$ Let $C$ be any
maximal chain containing $\left\{ z,y,w\right\} $. Then $C$ must also avoid $%
A,$ contradicting the assumption that $A$ is an antichain cutset. \ $\square 
$

\bigskip

\bigskip

\textbf{References}

\bigskip

[G] P.A. Grillet, Maximal chains and antichains, Fundamenta Math. 15 (1969)
157-167

\bigskip

[RZ] I. Rival, N. Zaguia, Antichain cutsets, Order 1 (1985) 235-247

\end{document}